\documentclass[draft]{publmathdeb}
\usepackage{amsmath,amsfonts,amssymb}
\newtheorem{definition}{Definition}
\newtheorem{example}{Example}
\newtheorem{theorem}{Theorem}

\begin{document}
\title{Geometry of Distributions \\ and $F-$Gordon equation}
\author{Mehdi Nadjafikhah}
\address{Iran University of Science and Technology,\\ Faculty
of Mathematics,\\ School of Mathematics.}
\email{m\_nadjafikhah@iust.ac.ir}
\urladdr{http://webpages.iust.ac.ir/m\_nadjafikhah}
\author{Reza Aghayan}
\address{Department of Mathematics,\\ Faculty of
Science,\\ Azad University, Eslamshahr Branch.}
\email{raghayan@yahoo.com}
%
%
\begin{abstract}
\noindent In this paper we describe the geometry of distributions
by their symmetries, and present a simplified proof of the
Frobenius theorem and some related corollaries. Then, we study the
geometry of solutions of $F-$Gordon equation; A PDE which appears
in differential geometry and relativistic field theory.
\end{abstract}
\submitted{December 28, 2008} \maketitle
\def\di{\displaystyle}
\def\rank{{\rm rank}\,}
\def\span{{\rm span}\,}
\def\Ann{{\rm Ann}\,}
\def\mod{{\rm mod}}
\def\Fl{{\rm Fl}}
\def\Ch{{\rm Char}}
\def\Sh{{\rm Shuf}}
\def\Sym{{\rm Sym}}
\def\RR{{\Bbb R}}
\def\cal{\bf}
\section{Introduction}
We begin this paper with the geometry of distributions. The main
idea here is the various notions of symmetry and their use in
solving a given differential equation. In section 2, we introduce
the basic notions and definitions.

In section 3, we describe the relation between differential
equations and distributions. In section 4, we present the geometry
of distributions by their symmetries, and find out the symmetries
of $F-$Gordon equation by this machinery. In section 5, we
introduce a simplified proof of the Frobenius theorem and some
related corollaries. In section 6, we describe the relations
between symmetries and solutions of a distribution.

In all steps, we study the $F-$Gordon equation as an application;
a partial differential equation which appears in differential
geometry and relativistic field theory. It is a generalized form
the Klein-Gordon equation "$u_{tt}-u_{xx}+u=0$". It is a
relativistic version of Schrodinger equation, which is used to
describe spinless particles. It was named after Walter Gordon and
Oskar Klein \cite{BEM, KRA}.
\section{Tangent and Cotangent Distribution}
Let $M$ be an $(m+n)-$dimensional smooth manifold.
\begin{definition}
A map ${\cal D}: M\rightarrow TM$ is called an {\it
$m-$dimensional tangent-distribution on $M$}, or briefly {\it
${\bf Tan}^m-$distribution}, if $${\cal D}_x:={\cal D}(x)\subseteq
T_xM \hspace{1cm} (x\in M)$$ is an $m-$dimensional subspace of
$T_xM$. The smoothness of $\cal D$ means that: \\ For each $x\in
M$, there exists an open neighborhood $U$ of $x$; and smooth
vector fields $X_1,\cdots, X_m$ such that:
\begin{eqnarray*} {\cal D}_y &=& \big<X_1(y), \cdots, X_m (y)\big> \\ &:=&
\span_\RR\{ X_1 (y), \cdots,X_m (y)\}\hspace{1cm} (y\in
U)\end{eqnarray*}
\end{definition}
\begin{definition}
A map $D:M\rightarrow T^{*}M$ is called an {\it $n-$dimension
cotangent-distribution on $M$}, or briefly {\it ${\bf
Cot}^n-$distribution}, if
$$D_x:=D(x)\subseteq T_x^{*}M\hspace{1cm} (x\in M)$$
is an $n-$dimensional subspace of $T^*_xM$. The smoothness of $D$
means that: \\ For each $x\in M$, there exists an open
neighborhood $U$ of $x$ and smooth $1-$forms $\omega^1,
\cdots\omega^n$ such that:
\begin{eqnarray*} D_y&=&\big<\omega^1(y), \cdots, \omega^n (y)\big>\\
&:=&\span_\RR\{ \omega^1 (y), \cdots,\omega^n(y)\} \hspace{1cm}
(y\in U) \end{eqnarray*}
\end{definition}
\paragraph{}
In the sequel, without loss of generality, we can assume these
definitions are globally satisfied.

\medskip There is a correspondence between these two types of
distributions. For ${\bf Tan}^m-$distribution $\cal D$, there
exist nowhere zero smooth vector fields $X_1,\cdots,X_m$ on $M$
such that: ${\cal D}=\big<X_1,\cdots,X_m\big>$; and similarly, for
${\bf Cot}^n-$distribution $D$, there exist global smooth
$1-$forms $\omega^1,\cdots,\omega^n$ on $M$ such that
$D=\big<\omega^1,\cdots,\omega^n\big>$.
\begin{example} \label{Exa:1} {\bf (Cartan distribution)} $\;\;$ \rm Let $M=\RR^{k+1}$. Denote the coordinates
in $M$ by $x,p_0,p_1 ,... ,p_k$ and given a function
$f(x,p_0,\cdots,p_{k-1})$ consider the following differential
$1-$forms:
\begin{eqnarray*} &&\omega^0=dp_0-p_1\,dx,\;\;\;
\omega^2=dp_1-p_2\,dx,\;\;\cdots
\\ &&\omega^{k-2}=dp_{k-2}-p_{k-1}\,dx,\;\; \omega^{k-1}=dp_{k-1}-
f(x,p_0,\cdots,p_{k-1})\,dx, \end{eqnarray*}
and the distribution $D=\big<\omega^0,\cdots,\omega^{k-1}\big>$.
This is the $1-$dimensional distribution, called the {\it Cartan
distribution}. This distribution can also be described by a single
vector field $X$, ${\cal D}=\big<X\big>$, where
$$X=\partial_x+p_1\,\partial_{p_0}+p_2\,\partial_{p_1}+\cdots+p_{k-1}\,\partial_{p_{k-2}}+
f(x,p_0,\cdots,p_{k-1})\,\partial_{p_{k-1}}.$$
\end{example}
\begin{example} \label{Exa:2} ({\bf $F-$Gordon equation}) \rm Let $F:\RR^5\to\RR$ be a
differentiable function. The corresponding $F-$Gordon PDE is
$u_{xy}=F(x,y,u,u_x,u_y)$. We construct $7-$dimensional
sub-manifold $M$ defined by $s=F(x,y,u,p,q)$, of
$$J^2(\RR^2,\RR)=\{{x,y,u,p=u_x,q=u_y,r=u_{xx},s=u_{xy},t=u_{yy}}\}.$$
Consider the $1-$forms
\begin{eqnarray*} \omega^1=du-p\,dx-q\,dy,\;\;\; \omega^2=dp-r\,dx-F\,dy, \;\;\;
\omega^3=dq-F\,dx-t\,dy. \end{eqnarray*}
This distribution can also be described by the following vector
fields
\begin{eqnarray*} X_1 =
\partial_x+p\,\partial_u+r\,\partial_p+F\,\partial_q,\;\;
X_2 =
\partial_y+q\,\partial_u+F\,\partial_p+t\,\partial_q,\;\;
X_3 = \partial_r,\;\; X_4 =
\partial_t.
\end{eqnarray*}
\end{example}
\begin{definition}
Let ${\cal D}:M\rightarrow TM$ be a ${\bf Tan}^m-$distribution and
set
$$\Ann{\cal D}_x :=\{\omega_x\in
T_x^*M\,\big|\,\omega_x\big|_{{\cal D}_x}=0\}.$$ It is clear that
$\dim\Ann{\cal D}_x=n$. An $1-$form $\omega\in \Omega^1(M)$ {\it
annihilates} $\cal D$ on a subset $N\subset M$, if and only if
$\omega_x\in\Ann{\cal D}_x$ for all $x\in M$.

The set of all differential $1-$forms on $M$ which annihilates
$\cal D$, is called {\it annihilator of $\cal D$} and denoted by
$\Ann{\cal D}$.
\end{definition}
\medskip Therefore, for each ${\bf Tan}^m$-distribution
\begin{eqnarray*} {\cal D}:M\rightarrow TM, \;\;\;\; {\cal D}:x \mapsto {\cal
D}_x, \end{eqnarray*}
we can construct a ${\bf Cot}^n-$distribution
\begin{eqnarray*} D:M\rightarrow T^*M, \;\;\;\; D:x \mapsto D_x=\Ann{\cal
D}_x;\end{eqnarray*}
and viceversa; In the other words, for each ${\bf
Tan}^m-$distribution ${\cal D}=\big<X_1, \cdots, X_m\big>$, we can
construct a ${\bf Cot}^n-$distribution $D=\Ann{\cal
D}=\big<\omega^1, \cdots ,\omega^n\big>$, and viceversa.
\begin{theorem} \label{Pro:1}
a) $\cal D$ and its annihilator, are modules over
$C^{\infty}(M)$. \\
b) Let $X$ be a smooth vector field on $M$ and $\omega\in
\Ann{\cal D}$, then $${\cal L}_X\omega\equiv-\omega \circ {\cal
L}_X\;\;\;\mod\;\;\;{\cal D}$$
\end{theorem}
{\it Proof:}  (a) is clear, and for (b), if $Y$ belongs to $\cal
D$, then $\omega(Y)=0$, and \\[2mm] \mbox{ } \hfill\ $({\cal
L}_X\omega)Y = X.(\omega(Y))-\omega[X.Y] = -\omega[X.Y]
=-(\omega\circ {\cal L}_X)Y.$ \hfill\ $\mbox{ }$
\section{Integral Manifolds, Maximal Integral Manifolds.}
\begin{definition}
Let $\cal D$ be a distribution. A bijective immersed sub-manifold
$N\subset M$, is called an {\it integral manifold of $\cal D$} if
one of the following equivalence conditions is satisfied:
\begin{itemize}
\item[1)]$T_x N\subseteq {\cal D}_x$, for all $x\in N$.
\item[2)]$N\subseteq \bigcap^{n}_{i=1}\ker \omega^i$.
\end{itemize}
$N\subset M$ is called {\it maximal integral manifold} if for each
$x\in N$, there exists an open neighborhood $U$ of $x$ such that
there is no integral manifold $N'$ containing $N\cap U$.
\end{definition}
\medskip It is clear that, the dimension of maximal integral
manifold does not exceed, the dimension of the distribution.
\begin{definition}
${\cal D}$ is called a {\it completely integrable distribution},
or briefly CID, if for all maximal integral manifold $N$, one of
the following equivalence conditions is satisfied:
\begin{itemize}
\item[1)] $\dim N= \dim {\cal D}$.
\item[2)] $T_xN={\cal D}_x$ for all $x\in N$
\item[3)] $N\subseteq\bigcap^n_{i=1}\ker\,\omega^i$, and if $N'$ be an
integral manifold with $N\cap N'\neq\emptyset$, then $N'\subseteq
N$.
\end{itemize}
\end{definition}
\paragraph{}
In the sequel, the set of all maximal integral manifolds is
denoted by $\cal N$.
\begin{theorem} ${\cal
N}=\bigcap^n_{i=1}\ker\,\omega^i$; that is $\omega^i|_{\cal N}=0$
for $i=1,\cdots ,n$.
\end{theorem}
\begin{example}{\bf (Continuation of Example \ref{Exa:1})} $\;\;$
\label{Exa:1c} \rm If $N$ is an integral curve of the distribution
then $x$ can be chosen as a coordinate on $N$, and therefore
$$N=\{(x,h_0(x),h_1(x),\cdots,h_{k-1}(x))\,|\,x\in\RR\}.$$
Conditions $\omega^0|_N=0,\cdots,\omega^{k-1}|_N=0$ imply that
$h_1=h_0', h_2=h_1',\cdots, h_{k-1}=h_{k-1}'$ or that
$$N=J^{k-1}h=\{(x,h(x),h'(x),\cdots, h^{(k-1)}(x)\,|\,x\in\RR\}$$ for
some function $h:\RR\to\RR$.

The last equation $\omega^{k-1}|_N=0$ gives us an ordinary
differential equation
$h^{(k)}(x)=f(x,h(x),h'(x),\cdots,h^{(k-1)}(x))$.

The existence theorem shows us once more that the integral curves
do exist, and therefore the Cartan distribution is a CID.
\end{example}
\begin{example}{\bf (Continuation of Example \ref{Exa:2})} $\;\;$
\label{Exa:2c} \rm This distribution in not CID, because there is
no $4-$dimensional integral manifold, and $\dim {\cal D}=4$. For,
if $N$ be a $4-$dimensinal integral manifold of the distribution,
then $(x,y,u,p)$ can be chosen as coordinates on $N$, and
therefore
\begin{eqnarray*} N\,:\, \Big\{ \, \begin{array}{l} q=h(x,y,u,p),\;r=l(x,y,u,p),\\
t=m(x,y,u,p),\; s=F(x,y,u,p,h). \end{array}\end{eqnarray*}
Condition $\omega^1|_N=0$ imply that $-p\,dx-h(x,y,u,p)\,dy+du=0$,
which is impossible.

By the same reason, we conclude that there is not any
$3-$dimensional integral manifold.

Now, if $N$ be a $2-$dimensinal integral manifold of the
distribution, then $(x,y)$ can be chosen as coordinates on $N$,
and therefore
\begin{eqnarray*}  N\,:\, \Big\{ \, \begin{array}{l}
u=h(x,y),\,p=l(x,y),\,q=m(x,y),\\
r=n(x,y),\,t=o(x,y),\, s=F(x,y,u,p,q). \end{array}\end{eqnarray*}
Conditions $\omega^1|_N=0$ and $\omega^2|_N=0$ imply that $l=h_x$,
$m=h_y$, $n=l_x=h_{xx}$, and $o=m_y=h_{yy}$.

The last equation $\omega^3|_N=0$ implies that
$h_{xy}=F(x,y,h,h_x,h_y)$. This distribution is not CID.
\end{example}
\section{Symmetries}
In this section, we consider a distribution ${\cal
D}=\big<X_1,\cdots ,X_m\big>=\big<\omega^1,\cdots ,\omega^n\big>$
on manifold $M^{n+m}$.
\begin{definition} A diffeomorphism $F:M\rightarrow M$ is called a {\it
symmetry of}
 $\cal D$ if $F_*{\cal D}_x={\cal D}_{F(x)}$ for all $x\in M$.
\end{definition}
\medskip Therefore,
\begin{theorem}\label{Pro:3}
{\it The following conditions are equivalent:
\begin{itemize}
\item[1)] $F$ is a symmetry of $\cal D$;
\item[2)] $F^{*}\omega^i$s determine the same distribution $\cal
D$; that is ${\cal
D}=\big<F^{*}\omega^1,\cdots,F^{*}\omega^n\big>$;
\item[3)] $F^{*}\omega^i\wedge\cdots\wedge \omega^n=0$ for
$i=1,\cdots,n$;
\item[4)] $F^{*}\omega^i=\sum^n_{j=1}\;a_{ij}\,\omega^j$, where $a_{ij}\in C^{\infty}(M)$;
\item[5)] $(F_*X_i|_x)\in {\cal D}_{F(x)}$ for all $x\in M$ and
$i=1,\cdots,n$; and
\item[6)] $F_*X_i= \sum^n_{j=1}b_{ij}\,X_j$, where $b_{ij}\in
C^{\infty}(M)$.
\end{itemize}}
\end{theorem}
\begin{theorem}
If $F$ be a symmetry of $\cal D$ and $N$ be an integral manifold,
then $F(N)$ is an integral manifold.
\end{theorem}
{\it Proof:} $F$ is a diffeomorphism, therefore $F(N)$ is a
sub-manifold of $M$. From other hand, if $x\in N$, then
$\omega^i|_{F(x)}=(F^{*}\omega^i)|_x=0$ for all $i=1,\cdots,n$,
therefore $F(N)=\{F(x)\,|\,x\in N\}$ is an integral manifold.
\hfill\ $\mbox{ }$
\begin{theorem}
Let $\cal N$ be the set of all maximal integral manifolds, and
$F:M\rightarrow M$ be a symmetry, then $F({\cal N})={\cal N}$.
\end{theorem}
{\it Proof:} If $x\in {\cal N}$, then
$\omega^i|_{F(x)}=(F^{*}\omega^i)|_x=0$ for all $i=1,\cdots,n$,
therefore $F(x)\in {\cal N}$; and $F({\cal N})\subset {\cal N}$.

Now if $y\in {\cal N}$, then there exists $x\in M$ such that
$F(x)=y$, since $F$ is a diffeomorphism. Therefore
$(F^*\omega^i)|_x=\omega^i|_{F(x)}=\omega^i|_y=0$ for all
$i=1,\cdots,n$; thus $x\in {\cal N}$ and ${\cal N}\subseteq
F({\cal N})$. \hfill\ $\mbox{ }$
\begin{definition}
 A vector field $X$ on $M$ is called an {\it infinitesimal symmetry of
distribution} $\cal D$, or briefly a {\it symmetry of ${\cal D}$},
if the flow $\Fl^X_t$ of $X$ be a symmetry of $\cal D$ for all
$t$.
\end{definition}
\begin{theorem}
{\it A vector field $X\in{\cal X}(M)$ is a symmetry if and only if
$${\cal L}_X\omega^i|_{\cal D}=0\;\;\;\mbox{for
all}\;\;\;i=1,\cdots,n.$$}
\end{theorem}
{\it Proof:} Let $X$ is a symmetry. If
$\Omega=\omega^1\wedge\cdots\wedge\omega^n$, then
$\{(\Fl^X)^*\omega^i\}\wedge \Omega=0$, by the (3) of Theorem
\ref{Pro:3}; Moreover by the definition ${\cal L}_X \omega^i :=
\frac{d}{dt}\big|_0 (\Fl^{X}_t)^*\omega^i$, one gets
\begin{eqnarray*} ({\cal L}_X \omega^i)\wedge \Omega &=&
\lim_{t\to0}\frac{1}{t}\,\big(
(\Fl^X_t)^*\omega^i-\omega^i\big)\wedge \Omega \\
&=& \di \lim_{t\rightarrow
0}\frac{1}{t}\big(\{(\Fl^X)^*\omega^i\}\wedge
\Omega-\omega^i\wedge \Omega^1\big) = 0. \end{eqnarray*}
Therefore ${\cal L}_X\omega^i|_{\cal {\cal D}}=0$.

In converse, let ${\cal L}_X\,\omega^i|_{\cal D}=0$ or ${\cal
L}_X\,\omega^i= \sum^n_{j=1}b_{ij}\;\omega^j$ for $i=1,\cdots,n$
and $b_{ij}\in C^\infty(M)$. Now, if
$\gamma_i(t):=\{(\Fl^X_t)^*\omega^i\}\wedge \Omega$, then
\begin{eqnarray}
\gamma_i(0)=\{(\Fl^X_0)^*\omega^i\}\wedge \Omega=0, \label{eq:1}
\end{eqnarray}
and
\begin{eqnarray*} \gamma'_i(t)&=&\frac{d}{dt}\{(\Fl^X_t)^*\omega^i\}\wedge
\Omega =((\Fl^X_t)^*{\cal L}_X\omega^i)\wedge \Omega\\
&=& (\Fl^X_t)^*\Big(\sum b_{ij}\,\omega^i\Big)\wedge \Omega = \sum
B_{ij}\,\{(\Fl^X_t)^*\omega^j\}\wedge \Omega \end{eqnarray*}
where $B_{ij}= (\Fl^X_t)^*b_{ij} = b_{ij}\circ \Fl^X_t$, and
\begin{eqnarray}
\gamma'_i(t)=\sum B_{ij}\,\gamma_i(t)\;\;,\;\;i=1,\cdots,n.
\label{eq:2}
\end{eqnarray}
Therefore,  $\gamma=(\gamma_1\cdots,\gamma_n)$ is a solution of
the linear homogeneous system of ODEs (\ref{eq:2}) with initial
conditions (\ref{eq:1}), and $\gamma$ must be identically zero.
\hfill\ $\mbox{ }$
\begin{theorem}
$X$ is symmetry if and only if for all $Y\in{\cal D}$ then
$[X,Y]\in{\cal D}$
\end{theorem}
{\it Proof:} By above theorem. $X$ is symmetry if and only if; for
all $\omega\in\Ann{\cal D}$ then ${\cal L}_X\omega\in\Ann{\cal
D}$.

The Theorem comes from the Theorem \ref{Pro:1}-(b). ${\cal
L}_X\omega=-\omega\circ {\cal L}_X$ on ${\cal D}$. In other words,
 $({\cal L}_X\omega)Y=-\omega[X,Y]$ for all $Y\in{\cal D}$.

Denote by $\Sym_{\cal D}$ the set of all symmetries of a
distribution ${\cal D}$.
\begin{example}{\bf (Continuation of Example \ref{Exa:1c})} $\;\;$
\label{Exa:1cc} \rm Let $k=2$. A vector field
$Y=a\,\partial_x+b\,\partial_{p_0}+c\,\partial_{p_1}$ is an
infinitesimal symmetry of $\cal D$ if and only if
$L_Y\omega^i\equiv0\,$mod$\,{\cal D}$, for $i=1,2$. These give two
equations
\begin{eqnarray*}  c = Xb-p_1\,Xa, \;\;\; Xc=f\,Xa+Yf, \end{eqnarray*}
\end{example}
\begin{example}{\bf (Continuation of Example \ref{Exa:2c})} $\;\;$
\label{Exa:2cc} \rm We consider the point infinitesimal
transformation:
\begin{eqnarray*} Z &=&
X(x,y,u)\,\partial_x+Y(x,y,u)\,\partial_y+U(x,y,u)\,\partial_u\\
&&
+P(x,y,u,p,q,r,t)\,\partial_p+Q(x,y,u,p,q,r,t)\,\partial_q\\
&& +R(x,y,u,p,q,r,t)\,\partial_r+T(x,y,u,p,q,r,t)\,\partial_t.
\end{eqnarray*}
Then, $Z$ is an infinitesimal symmetry of $\cal D$ if and only if
$L_Z\omega^i\equiv0\,$mod$\,{\cal D}$, for $i=1,2,3$. These give
ten equations
\begin{eqnarray*} && \hspace{-1cm} P_r=P_t=Q_r=Q_t=0,\\
&& \hspace{-1cm} p^2X_y+qpY_y+qU_x+pQ =pqX_x+qP+q^2Y_x+pU_y,\\
&& \hspace{-1cm} rpX_y+FY_y+qP_x+(qr-pF)P_p+(qF-pt)P_q+pXF_x+pYF_y+pUF_u\\
&& +pPF_p+pQF_q = qrX_x+qFY_x+pP_y+qR,\\
&& \hspace{-1cm} pFX_y+ptY_y+(qr-pF)Q_p+(qF-pt)Q_q+qQ_x+pT\\
&&=qFX_x+qtY_x+pQ_y+qPF_p+qYF_y+qXF_x+qQF_q+qUF_u,\\
&& \hspace{-1cm}Q_y+qQ_u+FQ_p+tQ_q=tY_y+tqY_u+FX_y+qFX_u+T,\\
&& \hspace{-1cm}U_y+qU_u=pX_y+pqX_u+qY_y+q^2Y_u+Q,\\
&& \hspace{-1cm}P_y+qP_u+FP_p+tP_q=rX_y+qrX_u+FY_y+qFY_u\\
&& +(XF_x+YF_y+UF_u+PF_p+QF_q).
 \end{eqnarray*}
Complicated computations using Maple, shows that:
\begin{eqnarray*} P&=&-pX_x-p^2X_u-qY_x-pqY_u+U_x+pU_u, \\
Q&=&\frac{1}{p}\Big(pqX_x-p^2X_y+q^2Y_x-pqY_y-qU_x+pU_y+qP\Big),\\
R &=& \frac{1}{q}\bigg(
\Big(pqX_x-p^2X_y+q^2Y_x-pqY_y-qU_x+pU_y+qP\Big).F_q\\
&&+((qr-pF).P_p
+(qF-pt).P_q+qP_x-pP_y\\
&&+pXF_x+pYF_y+pUF_u+pPF_p+(pY_y-qY_x).F-pr.(qX_x-pX_y) \bigg),
\\
T &=& \frac{1}{p^3}\bigg( (p^2t+q^2r-2pqF)P+p^2(pt+q^2F_q)X_x+p^2(qr-2pF-pqF_q)X_y\\
&&+q(q^2(r+pF_q)+3p(pt-qF))Y_x-p^2(q^2F_q+2pt+qF)Y_y\\
&&-(q^2r+pq^2F_q+p^2t+2pq)U_x+p^2qF_qU_y\\
&&-pq^2P_x+p^2qP_y+pq(pF-qr)P_p+pq(pt-qF)P_q\\
&&+pq(pXF_x+pYF_y+pUF_u+pPF_p+qPF_q)\\
&&-p^2q^2X_{xx}+2p^3qX_{xy}-p^4X_{yy}-pq^3Y_{xx}+2p^2q^2Y_{xy}-p^3qY_{yy} \\
&&+pq^2U_{xx}-2p^2qU_{xy}+p^3U_{yy} \bigg); \end{eqnarray*}
and $X=X(x,u-qy)$, $Y=Y(y,u-px)$ and  $U(x,y,u)$ must satisfy in
PDE:
\begin{eqnarray*} && \hspace{-5mm} (pF_p-F)X_x
+p(pF_p-2F)X_u+(qF_q-F)Y_y+q(qF_q-2F)Y_u\\
&&
-F_pU_x-F_qU_y+(F-pF_p-qF_q)U_u+U_{xy}+qU_{xu}+pU_{yu}+pqU_{uu}\\
&&=XF_x+YF_y+UF_u. \end{eqnarray*}
\end{example}
\section{A proof of Frobenius Theorem}
\begin{theorem}
 Let $X\in \Sym_{\cal D}\cap{\cal D}$, and $N$ maximal integral
manifold. Then $X$ is tangent to $N$.
\end{theorem}
{\it Proof:} Let $X(x)\not\in\;T_xN$. then exists open set $U$ of
$x$ and sufficiently small $\varepsilon$. such that
$\bar{N}:=\bigcup_{-\varepsilon<t<\varepsilon}\,\Fl^X_t\big(N\cap
U\big)$ is a smooth sub manifold of $M$.

Since $X\in{\cal D}$, So $\bar{N}$ is an integral manifold.

Since $X\in\Sym_{\cal D}$, So tangent to $\Fl^X_t\big(N\cap
U\big)$ belongs to ${\cal D}$, for all
$-\varepsilon<t<\varepsilon$.

On the other hand, tangent spaces to $\bar{N}$ are sums of tangent
spaces to $\Fl^X_t(N\cap U)$ and the $1-$dimensional subspace
generated by $X$, but both of them belongs to ${\cal D}$. and
there means $\bar{N}\subset N$
\begin{theorem}
If $X\in{\cal D}\cap\Sym_{\cal D}$ and $N$ be a maximal integral
manifold, then $\Fl^X_t(N)=N$ for all $t$.
\end{theorem}
\begin{theorem}[Frobenious]
A distribution ${\cal D}$ is completely integrable, if and only if
it is closed under Lie bracket. In other words, $[X,Y]\in{\cal D}$
for each $X,Y\in{\cal D}$.
\end{theorem}
{\it Proof:} Let $N$ be a maximal integral manifold with
$T_xN={\cal D}_x$. Therefore, for all $X,Y\in{\cal D}$ there $X$
and $Y$ are tangent to $N$ and so $[X,Y]$, is also tangent to $N$.

On the other hand, Let for all $X,Y\in{\cal D}$ there
$[X,Y]\in{\cal D}$. By Theorem, all $x\in{\cal D}$ is symmetry
too. and so all $X\in{\cal D}$ is tangent to $N$. and this means
$T_x\,N={\cal D}_x$, for all $x\in N$.
\begin{theorem}
A distribution ${\cal D}$ is completely integrable if and only if
${\cal D}\subset \Sym_{\cal D}$.
\end{theorem}
\begin{theorem}
Let ${\cal D}=\big<\omega^1,\cdots,\omega^n\big>$ be completely
Integrable distribution and $X\in{\cal D}$. Then the differential
$1-$forms $(\Fl_t^X)^*\omega^1,\cdots,(\Fl_t^X)^*\omega^n$ vanish
on ${\cal D}$ for all $t$.
\end{theorem}
{\it Proof:} ${\cal D}$ is completely integrable, then $X$ is
symmetry. Hence $$(\Fl_t^X)^*\omega^i=\sum_j\;a_{ij}\,\omega^j.$$
\section{Symmetries and Solutions}
\begin{definition}
If an (infinitesimal) symmetry $X$ belongs to the distribution
${\cal D}$, then it is called a {\it characteristic symmetry}.
Denote by $\Ch({\cal D}):=S_{\cal D}\cap {\cal D}$ the set of all
characteristic symmetries.
\end{definition}

\medskip It is shown that $\Ch({\cal D})$ is an ideal of the Lie
algebra $S_{\cal D}$, and is a module on ${\cal C}^\infty(M)$.
Thus, we can define the quotient Lie algebra
\begin{eqnarray*} \Sh({\cal D}) :=\Sym_{\cal D}/\Ch({\cal D}). \end{eqnarray*}
\begin{definition}
Elements of $\Sh({\cal D})$ are called {\it shuffling symmetries
of $\cal D$}.
\end{definition}

\medskip Any symmetry $X\in \Sym_{\cal D}$ generates a
flow on $\cal N$ (the set of all maximal integral manifolds of
$\cal D$), and, in fact the characteristic symmetries generate
trivial flows. In other words, classes $X\;{\rm mod}\;\Ch({\cal
D})$ mix or "shuffle" the set of all maximal manifolds.
\begin{example}{\bf (Continuation of Example \ref{Exa:1cc})} $\;\;$
\label{Exa:1ccc} \rm Let $k=2$. In this case $$\partial_x
\,\equiv\,-p_1\,\partial_{p_0}-f\,\partial_{p_1}\,\,{\rm
mod}\,\,\Ch({\cal D}),$$ therefore, $\Sh({\cal D})$ spanned by
$Z=(b-ap_1)\,\partial_{p_0}+(c-af)\,\partial_{p_1}$, where
\begin{eqnarray*} c = Xb-p_1\,Xa, \;\;\; Xc=f\,Xa+Yf.
 \end{eqnarray*}
\end{example}
\begin{example}{\bf (Continuation of Example \ref{Exa:2cc})} $\;\;$
\label{Exa:2ccc} \rm In this case, we have
\begin{eqnarray*}
& \partial_x  \equiv  -p\,\partial_u-r\,\partial_p-F\,\partial_q,\;\;\;\; \partial_r  \equiv  0, & \\
& \partial_y  \equiv  -q\,\partial_u-F\,\partial_p-t\,\partial_q,
\;\;\;\;\; \partial_t  \equiv  0, & \end{eqnarray*}
in $\Sh({\cal D})$. Therefore $\Sh({\cal D})$ spanned by
\begin{eqnarray*} W =
(U-pX-qY)\,\partial_u+(P-rX-FY)\,\partial_p+(Q-FX-tY)\,\partial_q,
\end{eqnarray*}
where
\begin{eqnarray*} P&=&-pX_x-p^2X_u-qY_x-pqY_u+U_x+pU_u, \\
Q&=&\frac{1}{p}\Big(pqX_x-p^2X_y+q^2Y_x-pqY_y-qU_x+pU_y+qP\Big),
\end{eqnarray*}
and $X=X(x,u-qy)$, $Y=Y(y,u-px)$ and  $U(x,y,u)$ must satisfy in
PDE:
\begin{eqnarray} && \hspace{-5mm} (pF_p-F)X_x
+p(pF_p-2F)X_u+(qF_q-F)Y_y+q(qF_q-2F)Y_u,\nonumber \\
&&
-F_pU_x-F_qU_y+(F-pF_p-qF_q)U_u+U_{xy}+qU_{xu}+pU_{yu}\label{eq:2-1}\\
&& +pqU_{uu}=XF_x+YF_y+UF_u. \nonumber \end{eqnarray}
\end{example}
\begin{example} {\bf (Quasilinear Klein-Gordon Equation)} \rm
In this example, we find the shuffling symmetries of Quasilinear
Klein-Gordon Equation:
\begin{eqnarray*} u_{tt}-\alpha^2\,u_{xx}+\gamma^2\,u=\beta\,u^3, \end{eqnarray*}
, as an application of previous example, where $\alpha$, $\beta $,
and $\gamma$ are real constants. The equation can be transformed
by defining $\xi=\frac{1}{2}(x-\alpha t)$ and
$\eta=\frac{1}{2}(x+\alpha t)$. Then, by the chain rule, we obtain
$\alpha^2\,u_{\xi\eta}+\gamma^2\,u=\beta\,u^3$. This equation
reduce to
\begin{eqnarray} u_{xy}=au+bu^3, \label{eq:3}\end{eqnarray}
by $t=y$, $a=-(\gamma/\alpha)^2$ and $b=\beta/\alpha^2$.

By solving the PDE (\ref{eq:2-1}), we conclude that $\Sh({\cal
D})$ spanned by the three following vector fields:
\begin{eqnarray*} X_1 &=& (px-qy)\,\partial_u - (p+yu^2(a+bu)-rx)\,\partial_p +
(q+xu^2(a+bu)-ty)\,\partial_q,\\
X_2 &=& q\,\partial_u + u^2(a+bu)\,\partial_p + t\,\partial_q,\\
X_3 &=& p\,\partial_u + r\,\partial_p + u^2(a+bu)\,\partial_q.
 \end{eqnarray*}
For example, we have
\begin{eqnarray*} && \hspace{-5mm}
\Fl_s^{X_3}(x,y,u,p,q,r,t)=\Big(x,y,u+sp+\frac{s^2}{2}r,p+sr,q+s.u^2(a+bu)\\
&&+\frac{s^2}{40}.up(2a+3bu)+\frac{s^3}{42}.(14ap^2+42bup^2+14aur+21bu^2r)\\
&&+\frac{s^4}{4}.p(bp^2+ar+3bur)
 +\frac{s^5}{20}.r(6bp^2+ar+3bur)+\frac{s^6}{8}.bpr^2+\frac{s^7}{56}.br^3,r,t\Big) \end{eqnarray*}
and if $u=h(x,y)$ be a solution of (\ref{eq:3}), then
$\Fl_s^{X_3}(x,y,h,h_x,h_y,h_{xx},h_{yy})$ is also a new solution
of (\ref{eq:3}), for sufficiently small $s\in{\bf R}$.
\end{example}

\label{lastpage}
\end{document}